\documentclass{article}
\usepackage{authblk}
\usepackage[scale=0.77]{geometry}
\usepackage{graphicx,xcolor} 
\usepackage[bf,font=small]{caption}
\usepackage{comment}
\usepackage{libertinus}
\usepackage{url}
\usepackage{amsfonts,amssymb,amsmath,amsthm,mathtools,thmtools, thm-restate, mathrsfs}

\declaretheorem[name=Theorem]{theorem}
\declaretheorem[name=Lemma, sibling=theorem]{lemma}

\newcommand{\eps}{\varepsilon}

\title{\textsc{On the Modular Chromatic Index\\ of Random Hypergraphs}}
\author[1]{Gaia Carenini}
\author[2]{Samuel Coulomb}
\affil[1]{DPMMS, University of Cambridge, Cambridge, United Kingdom.}
\affil[2]{Department of Computer Science, École Normale Supérieure -- PSL, Paris, France.}
\date{\today}

\begin{document}

\maketitle

\begin{abstract}
Let $k,r \ge 2$ be two integers. We consider the problem of partitioning the hyperedge set of an $r$-uniform hypergraph $H$ into the minimum number $\chi_k'(H)$ of edge-disjoint subhypergraphs in which every vertex has either degree $0$ or degree congruent to $1$ modulo $k$. For a random hypergraph $H$ drawn from the binomial model $\mathbf{H}(n,p,r)$, with edge probability $p \in (C\log(n)/n,1)$ for a large enough constant $C>0$ independent of $n$ and satisfying $n^{r-1}p(1-p)\to\infty$ as $n\to\infty$, we show that asymptotically almost surely $\chi_k'(H) = k$ if $n$ is divisible by $\gcd(k,r)$, and $\max(k,r) \le \chi_k'(H) \le k+r+1$ otherwise. A key ingredient in our approach is a sufficient condition ensuring the existence of a $k$-factor—a $k$-regular spanning subhypergraph—within subhypergraphs of a random hypergraph from $\mathbf{H}(n,p,r)$, a result that may be of independent interest.

Our main result extends a theorem of Botler, Colucci, and Kohayakawa (2023), who proved an analogous statement for graphs, and provides a partial answer to a question posed by Goetze, Klute, Knauer, Parada,  Peña, and Ueckerdt (2025) regarding whether $\chi_2'(H)$ can be bounded by a constant for every hypergraph $H$. 
\end{abstract}


\section{Introduction}
Decomposing the edge set of a graph or hypergraph into substructures subject to specific constraints is a classical and deeply studied problem in combinatorics. Beyond the extensively explored cases of matchings, stars, paths, and cycles, a prominent line of research focuses on \emph{degree-constrained decompositions}, partitions of the edge set in which each vertex must satisfy a prescribed degree condition within every part of the decomposition. A notable example is \emph{odd decomposition}, in which each subgraph in the partition is an \emph{odd graph} -- that is, a graph where every vertex has an odd degree.

The concept was introduced by Pyber in 1991, who defined the \emph{odd chromatic index} of a graph $G$, denoted $\chi'_2(G)$, as the minimum number of edge-disjoint subgraphs into which the edge set can be partitioned such that each subgraph is an odd graph \cite{pyber1991covering}. Pyber proved that for every simple graph $G$, $\chi'_2(G)\leq 4$. Since then, several works have further investigated odd edge decompositions and their combinatorial properties \cite{atanasov2016odd, kano2018decomposition, luvzar2014odd, matrai2006covering, petruvsevski2018odd, petruvsevski2019coverability, petruvsevski2021odd, petruvsevski2023odd}. 

A natural extension of odd decompositions arises by imposing \emph{modular degree constraints}. Given an integer $k\ge 2$, one considers partitions of the edge set into subgraphs in which each vertex has degree congruent to $1$ modulo $k$. This framework generalizes the notion of odd graphs, which correspond to the case $k=2$. For a graph $G$, the \emph{modular $k$-chromatic index}, denoted $\chi'_k(G)$, is defined as the minimum number of edge-disjoint subgraphs into which the edge set can be partitioned so that each subgraph satisfies the modular degree condition above. Scott~\cite{scott1997graph} showed that every graph $G$ admits a modular $k$-edge coloring with at most $5 k^2 \log k$ colors and asked whether the number of colors could in fact be bounded linearly in $k$. Subsequent work by Botler, Colucci, and Kohayakawa~\cite{Modk-198} established that $\chi'_k(G)\le 198 k - 101$ for all graphs $G$, and conjectured that the optimal bound should be of the form $k+C$ for some absolute constant $C$. In support of their conjecture, Botler, Colucci, and Kohayakawa~\cite{modk-random} showed that for \(k \ge 2\), a graph \(G\) drawn from the binomial random graph model \(\mathbf{G}(n,p)\) with \(p \in \bigl(C n^{-1}\log n, 1\bigr)\), where \(C > 0\) depends only on \(k\), and satisfying $n(1-p) \to \infty$ as $n \to \infty$, has $\chi'_k(G) = k$ or $\chi'_k(G) = k+1$, depending on the parity of $n$, with probability tending to 1 as $n \to \infty$. Progress towards confirming their conjecture was made by Nweit and Yang \cite{nweit2024mod}, who improved the multiplicative constant to \(177\), and later by Berthe, Bonamy, Botler, Carenini, Colucci, Dumas, Ghasemi, and Mariano Viana Neto \cite{berthe2025}, who further refined the bound to $\chi'_k(G)\le 9 k + o(k)$. 

 To date, the modular $k$-chromatic index has been investigated only in the setting of graphs, even in the case $k=2$. In this work, we initiate the study of the modular $k$-chromatic index in hypergraphs. For an integer $k \ge 2$, we call a hypergraph $H$ a \emph{$1_k$-hypergraph} if every vertex of $H$ has degree either $0$ or congruent to $1$ modulo $k$. The \emph{modular $k$-chromatic index} of a hypergraph $H$, denoted by $\chi_k'(H)$, is the minimum number of edge-disjoint $1_k$-subhypergraphs whose union covers all edges of $H$.

We focus on determining $\chi_k'(H)$ for a typical hypergraph $H$. Formally, given two integers $n \ge 1$ and $r \ge 2$ and a real number $p \in [0,1]$, let $\mathbf{H}(n,p,r)$ denote the probability space of $r$-uniform hypergraphs on the vertex set $[n] := \{1,2,\dots,n\}$ where each $r$-subset of vertices is chosen as a hyperedge independently with probability $p$. We write $H_{n,p,r}$ for a random point sampled from this space. We say that a hypergraph property $\mathcal{P}$ holds \emph{asymptotically almost surely} (a.a.s.) if the probability that $H_{n,p,r}$ verifies $\mathcal{P}$ tends to 1 as $n \to \infty$.

Our main result is the following:

\begin{theorem}[Main Theorem]\label{thrm:main}
    Let $k\ge2$ and $r\ge2$ be fixed integers, and let $d$ denote $\gcd(k,r)$. There exists a positive constant $C$ such that for every integer $n \ge 1$ and all $p \in \left(C\log(n)/n,1\right)$ such that $n^{r-1}p(1-p) \rightarrow \infty$ as $n \rightarrow \infty$, a.a.s. the mod $k$ chromatic index of  $H_{n,p,r}$ satisfies:
    \begin{enumerate}
        \item[(i)] if $n \equiv 0 \pmod d$, then $\chi_k'(H_{n,p,r})=k$, and 
        \item[(ii)] if $n \not\equiv 0 \pmod d$, then $\max(k,r) \le \chi_k'(H_{n,p,r}) \le k+r+1$.
    \end{enumerate}
\end{theorem}

The proof of Theorem 1 follows the general approach of Botler, Colucci, and Kohayakawa, who established an analogous result for graphs \cite{modk-random}. However, several modifications are required to accommodate the hypergraph setting, where certain tools—such as Hall’s theorem—no longer apply. In particular, our main technical contribution lies in identifying a sufficient condition that ensures the existence of a $k$-factor, that is a $k$-regular spanning subhypergraph, within subhypergraphs of a random hypergraph drawn from $\mathbf{H}(n,p,r)$. In particular, we prove the following result: 

\begin{theorem}[Main Technical Theorem]\label{thrm:factor}
    Let $r \ge 2$ be an integer and $\gamma,\eps,\eta\in(0,1)$ be real numbers. There exists a positive constant $C$ such that for every integer $n\ge1$ and all $p \in (C\log(n)/n,1)$, a.a.s. any subhypergraph $H \subseteq H_{n,p,r}$ of order $m = |V(H)|$ satisfying:
    \begin{enumerate}
        \vspace{-1mm}
        \item[(i)] $m \ge \gamma n$,
        \vspace{-1mm}
        \item[(ii)] $m \equiv 0 \pmod r$, and
        \vspace{-1mm}
        \item[(iii)] $\delta_{r-1}(H) \ge (\frac 1 2 +\eps)pm$,
        \vspace{-1mm}
    \end{enumerate}
    contains a $k$-factor for every positive integer $k \le \eps(1-\eta)mp$.
\end{theorem}

The conditions of Theorem \ref{thrm:factor} are fairly weak, suggesting that the result may be applicable in other settings and is, therefore of independent interest. Indeed, the existence of $k$-factors in graphs and hypergraphs has a long and rich history; for a survey, see \cite{fiorini1977edge,fiorini1978bibliographic,chung1981recent,akiyama1985factors,gould1991updating,gould2003advances,plummer2007graph,akiyama2011factors}. In particular, several studies have focused on factors in random graphs and hypergraphs; see, for instance, \cite{erdHos1966existence,shamir1981factors,shamir1984large,frieze1995perfect, Cooper_Frieze_Molloy_Reed_1996,krivelevich1997triangle,haber2007fractional,johansson2008factors,burghart2024sharp}. 

To conclude, we observe that Theorem \ref{thrm:main} shows that, in a typical $r$-uniform hypergraph $H$, the value of $\chi_k'(H)$ can be bounded by a function that grows linearly with $k$ and $r$, and that does not depend on the size of $H$. It is natural to wonder whether a similar bound holds for all $r$-uniform hypergraphs.  This question for $k=2$ was recently posed by Goetze, Klute, Knauer, Parada,  Peña, and Ueckerdt. \cite{goetze2025strong}.

\paragraph{Organization of the paper.}
The remainder of the paper is organized as follows. In Section \ref{sec: proof-overview}, we provide an outline of proof of our main result. In Section \ref{sec:factor}, we present the proof of Theorem \ref{thrm:factor}. The proof of Theorem \ref{thrm:main} is instead contained in Section \ref{sec:main-theo}. 

\section{Proof Overview}\label{sec: proof-overview}

We now outline the main ideas of the proof of Theorem~\ref{thrm:main} and Theorem \ref{thrm:factor}.

\vspace{3mm}

First, we establish that the mod-$k$ chromatic index of $H_{n,p,r}$ is at least $k$ by showing the existence of a vertex whose degree is both nonzero and divisible by $k$. This follows from standard concentration inequalities. The remainder of the proof then splits naturally into two cases, depending on the value of $n$. This is because a $r$-uniform $1_k$-hypergraph have order divisible by $d=\gcd(k,r)$, hence when $n \not\equiv 0 \pmod d$, we can't find a spanning $1_k$-subhypergaph in $H_{n,p,r}$, which makes this case considerably harder. 

If $n \equiv 0 \pmod d$, we start, as in the proof of the analogous results for graphs (Theorem 1 in \cite{modk-random}), by partitioning the vertices according to their degree modulo $k$. We show that, after deleting a small number of hyperedges, every subhypergraph induced by one of the degree classes in the partition satisfies the hypothesis of Theorem \ref{thrm:factor}. We apply this theorem to all these subhypergraphs, and identify a $i$-factor in the subhypergraph induced by vertices of degree congruent to $i+1$ modulo $k$. After removing these factors, all vertices have degree congruent to 1 modulo $k$, so we are left with a spanning $1_k$-subhypergraph. Lastly, we show that the hyperedges removed --both those discarded at the very beginning before applying Theorem 2 and those forming the factors -- can be partitioned into $k-1$ sets of independent edges, each of which forms a $1_k$-subhypergraph. 

If $n \not\equiv 0 \pmod d$, every $1_k$-subhypergraph must avoid at least one vertex. As in a typical $H_{n,p,r}$, every $(r-1)$-subset of vertices is contained in a hyperedge, we need at least $r$ $1_k$-subhypergraphs to cover all edges. For the upper bound, we randomly split the edges in $r+1$ subhyperhaphs $H_0,\dots,H_r$, each of order divisible by $d$. We identify, as in the previous case, a spanning $1_k$-hypergraph in $H_1,\dots,H_r$, and group all remaining hyperedges in $H_0$. Lastly, we partition $H_0$ in $1_k$-subhypergraphs like we did in the first case.

\vspace{3mm}

We can now turn our attention to Theorem \ref{thrm:factor}. Let $H\subseteq H_{n,p,r}$ be a subgraph verifying the hypothesis of Theorem \ref{thrm:factor}. In order to prove that $H$ contains a $i$-factor, we first show that it contains a perfect matching $M$. Then, it suffices to remove $M$ and iterate, as the union of $i$ disjoint perfect matchings is a $i$-factor. To this end, we consider an auxiliary bipartite graph $B$ in which a perfect matching correspond to a perfect matching in $H$. We partition at random $V(H)$ into $r$ equally-sized subsets $V_1, \dots, V_r$, ensuring that for all $(r-1)$-subsets of vertices $W$, the edges containing $W$ are well balanced between the parts. We randomly sample a labeling $\pi_i(1), \dots, \pi_i(m/r)$ of the vertices in each $V_i$. We contract the vertices $\pi_1(i), \dots, \pi_{r-1}(i)$ in a new vertex $i$, and let $B$ be the bipartite graph on $[m/r] \uplus V_r$ with an edge between $v \in V_r$ and $i \in [m/k]$ if $\{\pi_1(i), \dots, \pi_{r-1}(i),v\}$ is an hyperedge of $H$. We show that $B$ satisfies a strong expansion property, from which we deduce that $B$ verifies the hypothesis of Hall's Theorem (Theorem \ref{thrm:Hall}).

\section{Factors in Subhypergraphs of Random Hypergraphs}\label{sec:factor}

In this section, we prove Theorem \ref{thrm:factor}. We assume the reader has a basic familiarity with probability theory and make use of standard notations for Bernoulli and binomial random variables, which we briefly recall here. Let $n$ be a positive integer and $p \in [0,1]$. We write $X \sim \mathrm{Bern}(p)$ to denote a Bernoulli random variable with parameter $p$, and $X \sim \mathrm{Bin}(n,p)$ to denote a binomial random variable with parameters $n$ and $p$.

Given three integers $N, K, n$ with $0 \le n,K \le N$. We say that a random variable $X$ follows the \emph{hypergeometric distribution} with parameters $(N, K, n)$, denoted $X \sim \mathrm{Hypergeom}(N,K,n)$, if $\mathbb P(X = k) = \binom{K}{k} \binom{N-K}{n-k} / \binom{N}{n}$ for every non-negative integer $k$ such that $k \le K,n$ and $n+K-k \le N$. Roughly speaking, this distribution describes the number of elements from the set $[K]$ in a set $S$ obtained by sampling uniformly at random $n$ distinct elements in the set $[N]$.

To study how these random variables concentrate around their expected values, we employ the following version of the Chernoff bound.

\begin{lemma}[Chernoff Bound]\label{lem:Chernoff}
    Let $X$ be a random variable and $\mu = \mathbb E[X]$. If $X$ is the sum of independent Bernoulli random variables, or $X$ follows a binomial or hypergeometric distribution, then for all $a \in (0,1)$:\vspace{-1mm}
    \begin{itemize}
        \item[(i)]\label{lem:Chernoff_upper} \textbf{Upper Tail:}
        $\qquad \mathbb P(X \geq (1+a)\mu) \leq \exp \left( -a^2\mu/3 \right)$;
        \vspace{-1mm}
        \item[(ii)]\label{lem:Chernoff_lower} \textbf{Lower Tail:}
        $\qquad \mathbb P(X \leq (1-a)\mu) \leq \exp \left( -a^2\mu/2 \right)$.
    \end{itemize}
     \vspace{-2mm}
\end{lemma}

The main concentration tool in our proof is the following result due to McDiarmid \cite{McDIARMID_2002}. 

\begin{theorem}[McDiarmid Bound]\label{trhm:McDiarmid}
    Given a set $X$, we denote by $\mathfrak S(X)$ the set of all permutations of $X$. 
    Let $B_1, \dots, B_r$ be a family of finite sets and let $\Omega$ denote the cartesian product $\prod_{i=1}^r \mathfrak S(B_i)$. Let $c$ and $s$ be positive constants, and $h$ be a function from $\Omega$ to $\mathbb R_{\ge0}$ which for every $\pi = (\pi_1, \dots, \pi_r) \in \Omega$ satisfies the following conditions:
    \begin{enumerate}
        \item[(i)] for each $i \in [r]$, swapping two elements in $\pi_i$ can change the value of $h(\pi)$ by at most $c$, that is, for any transposition $\tau \in \mathfrak S(B_i)$ and $\pi' = (\pi_1, \dots, \tau \circ \pi_i, \dots, \pi_r)$, we have $|h(\pi)-h(\pi')| \le c$;
        \item[(ii)] there exist subsets $A_1, \dots, A_r$ of $B_1, \dots, B_r$ respectively, verifying $\sum_{i=1}^r |A_i| \le s \cdot h(\pi)$ and such that for any $\pi' \in \Omega$, if $\pi_i(a) = \pi_i'(a)$ for all $i \in [r]$ and $a \in A_i$, then $h(\pi') \ge h(\pi)$.
    \end{enumerate}
    Then, there exists a positive constant $C=C(c,s)$ such that, for $\pi$ sampled uniformly at random in $\Omega$, and any real number $t \ge 0$: 
    $$
    \mathbb P \left( h(\pi) \le \mu-t-C\sqrt\mu \right) \le 2 \exp \left( -t^2/C\mu \right) \qquad \text{where } \mu = \mathbb E[h(\pi)].
    $$
\end{theorem}
Another crucial component of our argument is Hall’s theorem, which we now recall.
\begin{theorem}[Hall's Theorem]\label{thrm:Hall}
    Let $G=(A\cup B,E)$ be a bipartite graph with $|A|=|B|=n$. Then, $G$ contains a perfect matching if and only if the following holds:
    \begin{enumerate}
        \item For all $X\subseteq A$ of size $|X|\leq n/2$ we have $|N(X)|\geq |X|$, and
        \item For all $Y\subseteq B$ of size $|Y|\leq n/2$ we have $|N(Y)|\geq |Y|$.
    \end{enumerate}
\end{theorem}

To present our proof precisely, we first introduce the following definitions. Let $G$ be a graph and $X,Y \subseteq V(G)$ be subsets of vertices, we denote by $e_G(X,Y)$ the number of edges between $X$ and $Y$ in $G$. For any vertex $v \in V(G)$, we define the neighbourhood of $v$, denoted $N_G(v)$, as the subset of vertices adjacent to $v$, that is, $N_G(v) = \{ u \in V(G) : uv \in E(G) \}$. Let $H$ be a hypergraph and $W,U \subseteq V(H)$ be subsets of vertices. The \emph{degree} of $W$, denoted $\deg_H(W)$, is the number of hyperedges $e \in E(H)$ that fully contain $W$, that is, those satisfying $W \subseteq e$. For any integer $i \ge 1$, we define the \emph{minimum $i$-degree} of $H$, denoted $\delta_i(H)$, as the minimum degree among all $i$-subsets of $V(H)$. Further, we denote by $\deg_H(W,U)$ the number of hyperedges $e \in E(H)$ such that $W \subseteq e$ and $e-W \subseteq U$.   

\begin{proof}[Proof of Theorem \ref{thrm:factor}]
    We start by partitioning uniformly at random the set $V(H)$ into $r$ equally-sized subsets $V_1, \dots, V_r$. This is always possible since by hypothesis $|V(H)| \equiv 0 \pmod r$. Note that for every $(r-1)$-subset $W \subseteq V(H)$ and index $i \in [r]$, the value of $\deg_H(W,V_i)$ is a random variable that depends on the partition of $V(H)$ sampled before. Additionally, $\deg_H(W,V_i)\sim \mathrm{Hypergeom}(m, \deg_H(W), m/r)$. In particular, its expected value is $\mu_0 = \mathbb E[\deg_H(W,V_i)] = \deg_H(W) / r$, and on the other hand, by hypothesis $\deg_H(W) \ge \delta_{r-1}(H) \ge (\frac 1 2+\eps) pm$. Hence, the Chernoff bound with parameter $a_0 = \frac{\eps}{1+2\eps}$, together with our assumption on $p$, implies the following chain of inequalities:
    $$
    \mathbb P \left( \deg_H(W,V_i) \le (1-a_0)  \mu_0 \right) 
    \le \exp \left( -\frac{a_0^2\mu_0}{2} \right)
    \le \exp \left( -\frac{a_0^2(\frac 1 2+\eps) \gamma C}{2r} \log(n) \right)
    \le n^{-Ca_0\eps\gamma/4r}.
    $$
    Applying the union bound, we deduce that the following upper bound holds for $C > 4r^2/a_0\eps\gamma$:
    $$
    \mathbb P \left( \exists W \in \binom{V(H)}{r-1} \text{ and } \exists i \in [r]: \deg(W,V_i) \le (1-a_0)\frac{\deg_H(W)}{r} \right)
    \le r \binom{n}{r-1}n^{-Ca_0\eps\gamma/4r} < 1.
    $$
    Consequently, there exists a partition $\{V_1,\dots, V_r\}$ that satisfies $\deg_H(W,V_i) \ge (1-a_0) \deg_H(W)/r \ge (\frac 1 2 + \frac \eps 2)pm/r$ for all $(r-1)$-subsets $W \subseteq V(H)$ and all $i \in [r]$. From now on, we fix such a partition.

    For any family of permutations $\pi=\{\pi_1, \dots, \pi_{r-1}\}$ where each $\pi_i$ maps $[m/r]$ to $V_i$, we define a bipartite graph $B_\pi=B_\pi(H)$ as follows. Let $V(B_\pi) = [m/r] \uplus V_r$ with an edge between $i \in [m/r]$ and $v \in V_r$ in $B_\pi$  if and only if $\{\pi_1(i), \dots, \pi_{r-1}(i), v\}\in E(H)$. Remark that we have selected the partition $\{V_1,\dots, V_r\}$ so that $\deg_{B_\pi}(i) = \deg_H(\{\pi_1(i), \dots, \pi_{r-1}(i)\}, V_r) \ge (\frac 1 2+\frac \eps 2) mp/r$ for every $i \in [m/r]$. Next, we show that there exists a family of permutation $\pi$ satisfying $\delta(B_\pi) \ge (\frac 1 2+\frac \eps 4) mp/r$. In fact, we prove a stronger statement: such a family of permutations can be obtained by a uniform random sampling. Let $\pi$ be a family of permutations sampled uniformly at random from the cartesian product $\Omega = \prod_{i=1}^{r-1} \mathfrak S(V_i)$. Given a vertex $v \in V_r$, let $h_v(\pi) = \deg_{B_\pi}(v)$ and $\mu_1 = \mathbb E[h_v(\pi)]$. Using linearity of expectation, we derive that:
    $$
    \mu_1 %= \sum_{i=1}^{m/r} \mathbb P \left( \{\pi_1(i), \dots, \pi_{r-1}(i), v\} \in E[B_\pi] \right)
    = \sum_{i=1}^{m/r} \mathbb E \left[ \frac{\deg_H(\{\pi_2(i), \dots, \pi_{r-1}(i), v\})}{|V_1|} \right] \ge \left( \frac 1 2 + \frac \eps 2 \right) \frac{mp}{r}.
    $$
    Note that for each $i \in [r]$, swapping two elements in $\pi_i$ only affect two vertices in $B_\pi$, thus it changes the value of $h_v(\pi)$ by at most 2. Moreover, for any family of permutations $\pi' \in \Omega$, if the equality $\pi_i(u) = \pi_i'(u)$ holds for all $i \in [r]$ and all $u \in N_{B_{\pi}}(v)$, then the vertices in $N_{B_\pi}(v)$ are still neighbors of $v$ in $B_{\pi'}$, so $h_v(\pi') \ge h_v(\pi)$. We may now apply McDiarmid's inequality with parameters $c=2$, $s=r$, and $t=a_1\mu_1$ where $a_1=\frac{\eps}{3(1+\eps)}$, obtaining that there exists a positive constant $C_1$ such that:
    $$
    \mathbb P \left( h_v(\pi) \le (1-a_1) \mu_1 - C_1\sqrt{\mu_1} \right) \le 2 \exp \left( -\frac{a_1^2\mu_1}{C_1} \right).
    $$
    For $n$ large enough, this implies:
    $$
    \mathbb P \left( h_v(\pi) \le (\tfrac 1 2+\tfrac \eps 4)\frac{mp}{r} \right) \le 2 \exp \left( -\frac{a_1^2mp}{2C_1r} \right) < \frac 1 m.
    $$
    Applying the union bound, we deduce the following upper bound:
    $$
    \mathbb P \left( \exists v\in V_r : \deg_{B_\pi}(v) \le \left( \tfrac 1 2+\tfrac \eps 4 \right) \frac{mp}{r} \right) < \frac{|V_r|}{m} < 1.
    $$
    Hence, there exists a family of permutations $\pi$ such that $\delta(B_\pi) \ge (\frac 1 2+\frac \eps 4) mp/r$. From now on, we fix such a family.
    
     Consider two subsets $X \subseteq [m/r]$ and $Y \subseteq V_r$ both of size $x \le m/2r$, and let $t=(\frac 1 2+\frac \eps 5)mpx/r$. Observe that if $e_{B_\pi}(X,Y) > t$, then there exist a collection of $x$ disjoint $(r-1)$-subsets of vertices, say $W_1,\dots,W_x \subseteq V(H)$, in which no subset intersects $Y$ and such that the number of hyperedges in $H_{n,p,r}$ of the form $W_i \cup \{a\}$ with $a \in Y$ is larger than $t$. We now prove that a.a.s. $H_{n,p,r}$ does not contain such a collection. Given a collection of disjoint $(r-1)$-subsets $W_1, \dots, W_x \subseteq [n]$, let $Z(H_{n,p,r})$ denote $\sum_{i=1}^x \deg_{H_{n,p,r}}(W_i,Y)$. We observe that $Z(H_{n,p,r}) \sim Bin(x^2,p)$, and consequently that $\mu_2 = \mathbb E[Z(H_{n,p,r})] = x^2p \le mpx/2r<t$.

     We distinguish three cases according to the value of $x$. We start by assuming that $x \le mp/2r$. In this case, we can immediately conclude since $Z(H_{n,p,r}) \le x^2 < mpx/2r < t$. We can now suppose that $mp/2r \le x \le m/4r$. If this is the case, we apply Stirling's formula to obtain the following bound:
    $$
    \mathbb P (Z(H_{n,p,r}) \ge t) \le \binom{x^2}{t}p^t \le \left( \frac{ex^2p}{t} \right)^t
    \le \exp \left( -\left(\tfrac 1 2+\tfrac \eps 5\right)\frac{mpx}{r} \right) \le \exp \left( -\frac{\eps^2mpx}{50r} \right).
    $$
    On the other hand, if $x \ge m/4r$, then $\mu_2\ge mpx/4r$. In this case, by the Chernoff bound, we get:
    $$
    \mathbb P (Z(H_{n,p,r}) \ge t ) = \mathbb P \left( Z(H_{n,p,r}) \ge \left( 1+\tfrac{2\eps}{5} \right) \mu_2 \right) \le \exp \left( -\frac{4\eps^2 \mu_2}{50} \right) \le \exp \left( -\frac{\eps^2mpx}{50r} \right).
    $$
    By union bound, we obtain that the following holds whenever $C>50(r+1)r/\eps^2\gamma$:
    \begin{align*}
        &\mathbb P \left( \exists W_1, \dots, W_x \in \binom{V(H_{n,p,r})}{r-1}, \exists Y \in \binom{V(H_{n,p,r})}{x} : \sum_{i=1}^x \deg_{H_{n,p,r}}(W_i,Y) \ge t \right) \\
        &\le \sum_{x=mp/2r}^{m/2r} \binom{n}{r-1}^x \binom{n}{x} \exp \left( -\frac{\eps^2mpx}{50r} \right) \\
        &\le \sum_{x=1}^{n} n^{(r-1)x}n^x\exp \left( -C\frac{\eps^2\gamma x}{50r} \log(n) \right) \\
        &\le \sum_{x=1}^{n} n^{x (r-C\eps^2\gamma/50r)} \le n^{1+r-C\eps^2\gamma/50r} < 1.
    \end{align*}
    We conclude that a.s.s. for all $X \subseteq [m/r]$ and $Y \subseteq V_r$ both of size $x \le m/2r$, we have $e_{B_\pi}(X,Y) \le (\frac 1 2+\frac \eps 5)mpx/r$.    
    
    We now show that $B_{\pi}$ satisfies the hypothesis of Hall's Theorem. 
    Suppose by contradiction that there exists a subset of vertices $X \subseteq [m/r]$ (resp. $X \subseteq V_r$) of size $x \le m/2r$ such that $|N_{B_\pi}(X)| < |X|$, and let $Y$ be an arbitrary subset of $V_r$ (resp. $[m/r]$) of size $x$ fully containing $\mathcal N_{B_\pi}(X)$. We derive that:
    $$\left(\frac 1 2+\frac \eps 5\right)\frac{mpx}{r} \ge e_{B_\pi}(X,Y) \ge \delta(B_\pi)|Y| \ge \left(\frac 1 2+\frac \eps 4\right)\frac{mpx}{r},$$
    a contradiction. We can therefore apply Hall's theorem to $B_\pi$ and deduce that it contains a perfect matching. It is easy to see that one can obtain from a perfect matching $\{\{i,v_i\} : i \in [m/r]\}$ in $B_\pi$, a perfect matching $\{\{\pi_1(i), \dots, \pi_{r-1}(i), v_i\} : i \in [m/r]\}$ in $H_{n,p,r}$. 
    
    So far, we have proved that there exists a positive constant $C$ such that a.a.s. any subhypergraph $H \subseteq H_{n,p,r}$ verifying conditions (i), (ii), and (iii) contains a perfect matching $M$. We can now build a $k$-factor by repetitively identifying perfect matchings as follows. Let $H_1=H$ and $M_1$ be a perfect matching in $H$ identified as above. We define $H_2$ to be $H_1 - M_1$, and observe that also $H_2$ satisfies the sufficient conditions for $\eps' = \eps\eta$, thus guarantying the existence of a perfect matching $M_2$ in $H_2$. Then $M_1\cup M_2$ is a $2$-factor of $H$. More generally, for every positive integer $k' < \eps(1-\eta)mp$ and $k'$-factor $F$, we can claim the existence of a perfect matching $M'$ in $H-F$, and $F \cup M'$ is a $(k'+1)$-factor of $H$. This concludes the proof of the theorem. 
\end{proof}

\section{Proof of the Main Theorem (Theorem \ref{thrm:main})}\label{sec:main-theo}
In this section, we prove our main result. For this proof, we need the following lemma that is a mild variant of several known results, and for which a simple and elegant proof can be found in \cite{modk-random}.  

\begin{lemma}\label{lem:bin_modk}
    Let $k\ge2$ be an integer. There exists a positive constant $C$ such that for any two integers $n\ge1$ and $t$, any real number $p \in (0,1)$, and any random variable $X \sim \mathrm{Bin}(n,p)$, we have:
    $$
    \left|\mathbb P(X \equiv t \hspace{-2mm} \pmod k) - \frac 1 k \right| \leq \exp(-Cnp(1-p)).
    $$
\end{lemma}

\begin{proof}[Proof of Theorem \ref{thrm:main}]
    Given a hypergraph $H$ and an integer $i \in [k]$, let $V_i(H) \subseteq V(H)$ be the subset of vertices of degree congruent to $i$ modulo $k$. We consider an arbitrary subset $U \subseteq [n]$ of size $n-m$ where $m=\left\lfloor \frac{n}{11} \right\rfloor$. We sample a hypergraph $H \sim \textbf{H}(n,p,r)$ in two steps. We first generate the hyperedges of $H[U]$, i.e., we sample the hyperedges within $\binom{U}{r}$, and fix them. For every vertex $u \in U$ and index $i\in[k]$, let $X_i^u$ be a random variable accounting for the indicator function $\mathbf 1[u \in V_i(H)]$. These variables are independent since we have fixed the hyperedges among the vertices of $U$, and each $X_i^u$ follows a Bernoulli distribution with parameter $p_i^u = \mathbb P(\deg_H(u) \equiv i \pmod k)$. On the other hand, the value of $\deg_H(u)-\deg_{H[U]}(u)$ follows the distribution $\mathrm{Bin}(\binom{m}{r-1},p)$. Therefore, by Lemma \ref{lem:bin_modk}, there exists a positive constant $C_0$ such that $|p_i^u-1/k|\le \exp(-C_0\binom{m}{r-1}p(1-p))$. Further, since by hypothesis $n^{r-1}p(1-p) \rightarrow \infty$, we have $|p_i^u-1/k| \le 1/12k$ for all sufficiently large values of $n$. We now consider the random variable $X_i = |V_i(H) \cap U|=\sum_{u\in U} X_i^u$. Let $\mu_i = \mathbb E[X] = \sum_{u \in U} p_i^u$, and observe that:
    $$
    \frac{5n}{6k} \le \frac{11}{12k}|U| \le \mu_i \le \frac{13}{12k}|U| \le \frac{n}{k}.
    $$
    By applying both Chernoff bounds, it follows that $\mathbb P(|X_i-\mu_i| \ge \mu_i/10) \le 2\exp(-\mu_i/300)$, which in turn implies $\mathbb P(|X_i-n/k| \ge n/4k) \le 2\exp(-n/360k)$. A simple union bound allows us to deduce that a.a.s., for all $i \in [k]$, we have $3n/4k \le |V_i(H_{n,p,r})| \le 5n/4k$.
    
    Consequently, a.a.s. there exists a vertex $v\in V_k(H)$. Moreover, the probability that a vertex has degree $0$ in $H$ is at most $n(1-p)^n \le n (1-C\log(n)/n)^n \le n^{1-C}$ which, for $C>1$, goes to 0 as $n$ goes to infinity. Hence, one need at least $k$ edge-disjoint $1_k$-hypergraphs to cover the hyperedges incident to $v$, so a.a.s. $\chi_k'(H) \ge k$.
  
    Let $W \subseteq [n]$ be an arbitrary $(r-1)$-subset and denote by $E^W$ the set of all $r$-subsets containing $W$, that is, $E^W = \{e \in \binom{[n]}{r} : W \subset e\}$. As before, we sample a hypergraph $H \sim \textbf{H}(n,p,r)$ in two steps. We first generate the hyperdeges $e \in \binom{[n]}{r} - E^W$, and fix them. Let $X_i^W = \deg_H(W,V_i(H))$. Clearly $X_i^W = \deg_H(W,V_{i-1}(H-E^W))$, and the random variable $X_i^W$ follows the binomial distribution $\mathrm{Bin}(|V_{i-1}(H-E^W)|,p)$. Let $\mu_i^W = \mathbb E[X_i^W] = |V_{i-1}(H-E^W)|p$. By applying the Chernoff bound, since $|V_{i-1}(H-E^W)| \ge 3/4k$, it follows that a.a.s.:
    $$
    \mathbb P \left( X_i^W \le \frac{11}{12} \mu_i^W \right) \le \exp \left( -\frac{\mu_i^W}{300} \right) \le \exp \left( -\frac{C\log(n)}{400k} \right)
    $$
    By union bound, for $C > 400kr$, we deduce:
    $$
    \mathbb P \left( \exists W \in \binom{[n]}{r-1} \text{ and } \exists i \in [k] : \deg_H(W,V_i(H)) \le \frac{11}{16}np \right) \le kn^{r-1}n^{-C/400k}< \frac k n.
    $$
    Consequently, a.a.s. for every $(r-1)$-subset of vertices $W \subseteq [n]$ and every index $i \in [k]$, it holds that: $$\deg_{H_{n,p,r}}(W,V_i(H_{n,p,r})) \ge \frac{11}{16k}np \ge C\frac{11}{16k} \log(n).$$ In particular, for any subset of vertices $A \subseteq [n]$ whose size is bounded by a constant independent of $n$, we can a.a.s. find a hyperedge of the form $W\cup v$ with $v\in V_i(H_{n,p,r}) - A$. The constructions carried out in the next paragraphs heavily rely on this fact.
    In what follows, we distinguish two cases according to the residue of $n$ modulo $d = \gcd(k, r)$. 
    
   \paragraph{Case 1: $n \equiv 0 \pmod d$.} The first case we consider is $n \equiv 0 \pmod d$, and corresponds to the setting in which one may identify a spanning color class. Indeed, every $r$-uniform $1_k$-hypergraph $\tilde H$ has to verify the following conditions: $\sum_{v \in V(\tilde H)} \deg_{\tilde H} (v) \equiv |V(\tilde H)| \pmod k$ and $\sum_{v \in V(\tilde H)} \deg_{\tilde H}(v) \equiv 0 \pmod r$; thus, $|V(\tilde H)| \equiv 0 \pmod d$. For each integer $i \in [3,k]$, let $q_i$ be the remainder of $|V_i(H)|$ modulo $r$. We construct a family of disjoint stars $S_{i,j}$ for all integers $i \in [3,k]$ and $j \in [1,q_i]$, such that each star $S_{i,j}$ has $i-1$ branches, its center in $V_i(H)$, and all its leaves are $V_2(H)$. Let $G$ be the union of all stars $S_{i,j}$, and $H' = H - E(G)$. Remark that for each vertex $v \in V(G)$, we have $\deg_{H'}(v) \equiv 1 \pmod k$, and for each integer $i \in [3,k]$, as $|V_i(H) \cap V(G)| = q_i$, we have  $|V_i(H')| \equiv 0 \pmod r$.
    
    For each $i \in \{1,2\}$, let $q_i$ be the remainder of $|V_i(H')|$ modulo $r$. On the one hand, we have $q_1+q_2 \equiv \sum_{i=1}^k |V_i(H')| \equiv n \pmod r$, and by assumption $n \equiv 0 \pmod d$, thus $q_1+q_2 \equiv 0 \pmod d$. On the other hand, we have $q_1+2q_2 \equiv \sum_{i=1}^k i|V_i(H')| \pmod r$, and $\sum_{i=1}^k i|V_i(H')| \equiv r |E(H')|  \pmod k$, hence $q_1+2q_2 \equiv 0 \pmod d$. This implies that $q_2 \equiv 0 \pmod d$, and that there exists two positive integers $a<r$ and $b<k$ such that $ak=br+q_2$. We fix an arbitrary $(r-1)$-subset $W' \subseteq V_2(H')$ and construct a family $S_2$ of $ak$ hyperedges of the form $W' \cup \{v\}$ with $v \in V_2(H')$. Let $H''= H' - E(S_2)$. Remark that, for every vertex $v \in W'$, we have $\deg_{H''}(v) \equiv 2 \pmod k$, and for every vertex $v \in V(S_2)-W'$, we have $\deg_{H''}(v) \equiv 1 \pmod k$. As $|V(S_2)-W'| = ak = br + q_2$, we have $|V_2(H'')| \equiv 0 \pmod r$. Let $C_0 = |V(G)|+|V(S_2)|$, and note that $C_0 \le k^2r^2$.
    
    We proved that, for every integer $i \in [2,k]$, a.a.s. the subhypergraph $H_i'' \subseteq H''$ induced  by $V_i(H'')$ verifies: 
    \begin{enumerate}
        \item[(i)] $|V(H_i'')| \ge |V_i(H)|-C_0 \ge \frac{2}{3k}n$, and
        \item[(ii)] $|V(H_i'')| \equiv 0 \pmod r$, and
        \item[(iii)] $\delta_{r-1}(H_i'') \ge \frac{11}{16k}np-C_0 \ge (\frac 1 2+\frac{1}{40})|V(H_i'')|p$.
    \end{enumerate}
    By Theorem \ref{thrm:factor}, $V_i(H'')$ contains a $(i-1)$-factor $F_i$ that is the union of $i-1$ perfect matching $M_{i,1},\dots,M_{i,i-1}$. Observe that $G$ and $S_2-W'$ are disjoint and included in $V_1(H'')$, thus they are also disjoint from all factors, whereas the vertices of $W'$ are in $V_2(H'')$, hence they are covered by the factor $F_2$.
    
    We now detail how to partition the hyperedges of $H_{n,p,r}$ into $1_k$-subhypergraphs $H_1,\dots,H_k$. Put each matching $M_{i,j}$ into $H_j$, and for each star $S_{i,j}$, put its $i-1$ branches in $H_1,\dots,H_{i-1}$ respectively. Observe that every vertex has degree 0 or 1 in each of $H_1,\dots,H_{k-1}$. Now, put $S_2$ in $H_1$. The vertices in $W'$ have degree $ak+1$ in $H_1$, and the other vertices in $V(S_2)$ have degree 1 in $H_1$. Lastly, put all remaining edges in $H_k$. Every vertex has degree congruent to 1 modulo $k$ in $H_k$. Henceforth, when $n \equiv 0 \pmod  d$, a.a.s. $\chi_k'(H_{n,p,r})=k$.

    \paragraph{Case 2: $n \not\equiv 0 \pmod d$.}
    We now tackle the case where $n \not\equiv 0 \pmod d$. Suppose the edges of $H_{n,p,r}$ can be partitioned into $1_k$-subhypergraphs $H_1, \dots, H_{r-1}$. Recall that a $1_k$-hypergraph $H'$ with no isolated vertices must have order $|V(H')|$ congruent to 0 modulo $d$. Hence, for each index $i \in [r-1]$, there exists a vertex $v_i$ that is isolated in $H_i$. Let $W \subseteq [n]$ be any $(r-1)$-subset containing $\{v_1, \dots, v_{r-1}\}$. A.a.s. $W$ has positive degree, but an hyperedge containing $W$ cannot belong to any $H_i$, a contradition. Therefore, a.a.s. $\chi_k'(H_{n,p,r}) \ge r$.
    
    Let $q$ be the remainder of $n$ modulo $d$. Fix $r+1$ disjoint subsets $R_0, \dots, R_r \subseteq [n]$ of size $q$, and let $R = \bigcup_{i=0}^r R_i$. We first split the hypergraph $H$ into $r+1$ subhypergraphs $H_0, \dots, H_r$ where each $H_i$ has vertex set $V(H_i) = [n] - R_i$. Observe that no $r$-set can intersect every $R_i$, so for every $e \in \binom{[n]}{r}$, there is at least one index $i \in [0,r]$ such that $e \subseteq V(H_i)$. For each hyperedge $e \in E(H)$, we sample uniformly at random an index $i \in [0,r]$ such that $e \subseteq V(H_i)$ and add $e$ to $H_i$. All these samples are independent. Clearly, an $r$-set $e \subseteq V(H_i)$ is a hyperedge of $H_i$ with probability $p/j$, independently of other $r$-sets, where $j \in [0,r]$ is the number of index $i' \in [0,r]$ such that $e \subseteq V(H_{i'})$. It easy to show that the degree classes $V_1(H_i), \dots, V_k(H_i)$ in $H_i$ satisfy the same properties that we proved in $H_{n,p,r}$. More precisely, for every integers $i,i' \in [0,r]$ and $j \in [k]$, we have $3n/4k \le |V_j(H_i)| \le 5n/4k$, and for all $(r-1)$-subset $W \subseteq [n]$, it holds that $\deg_{H_i}(W,V_j(H_{i'})) \ge \frac{11}{16k}np$.

    For each index $i \in [r]$ and integer $u \in R_0$, let $d_i^u \in [k]$ be the index such that $u \in V_{d_i^u}(H_i)$. We construct a star $S_i^u \subseteq H_i$ centered in $u$, with $d_i^u-1$ branches, and whose leaves are in $V_{r+2}(H_{n,p,r})$ and disjoint from other stars. Let $S_i$ denote $\bigcup_{u \in R_0} S_i^u$, and remark that the vertices of $R_0$ have degree congruent to 1 modulo $k$ in  $H_i-S_i$. As in Case 1, we construct a forest of stars $G_i$ disjoint from $S_i$ so that for all integers $j \in [2,k]$, it holds that $|V_j(H_i-S_i-G_i)| \equiv 0 \pmod r$, and use Theorem \ref{thrm:factor} to find a $(j-1)$-factor $F_{i,j}$ in $V_j(H_i-S_i-G_i)$. Let $F_i=\bigcup_{j \in [k]} F_{i,j}$, and denote $H_i'$ for $H_i-S_i^W-G_i-F_i$. Note that $H_i'$ is a $1_k$-hypergraph.
    
    Let $S = \bigcup_{i \in [k]} S_i$. Recall that for each $i \in [k]$, we have $R_0 \subseteq V_1(H_i -S_i)$, so $R_0$ is disjoint from $G_i$ and $F_i$. Thus, we may consider $H_0' = H_0 + \bigcup_{i \in [k]}(G_i+F_i)$. Observe that for any vertex $v \notin V(S) \cup R_0$  and integer $j \in [k]$, if $v$ has degree congruent to $j$ modulo $k$ in $H_0'$, then it has degree congruent to $j+r$ modulo $k$ in $H_{n,p,r}$, as $v$ has degree congruent to $1$ modulo $k$ in $H_1',\dots,H_r'$. Moreover, the size of $V(S) \cup R_0$ is bounded by $C_0=dkr$ which is independent of $n$. Thus, we have $3n/4k - C_0 \le |V_j(H_i)| \le 5n/4k + C_0$, and for all $(r-1)$-subset $W \subseteq [n]$, it holds that $\delta_{l-1}^{H_i}(W,V_j(H_{i'})) \ge \frac{11}{16k}np - C_0$. We again construct a forest of stars $G_0$ disjoint from all previous stars so that for all $j \in [k]$, we have $|V_j(H_0'-G_0)| \equiv 0 \pmod r$. We then apply Theorem \ref{thrm:factor} to find for each $j \in [k]$, a $(j-1)$-factor $F_{0,j}$ in $V_j(H_0'-G_0)$. The vertices of $V(S)-R_0$ were chosen in $V_{r+2}(H_{n,p,r})$, thus they belong to $V_1(H_0')$. Hence, $F_0$, $S$, and $G_0$ are disjoint, further, they can easily be partitioned into $k$ $1_k$-subhypergraphs. Lastly, note that $H_0'-S-F_0-G_0$ is a $1_k$-hypergraph. Henceforth, a.a.s. $\chi_k'(H_{n,p,r}) \le k+r+1$.

    This concludes the proof of the theorem. 
\end{proof}

\section*{Acknowledgments}
We are grateful to Meike Hatzel and Marcelo Garlet Milani for bringing to our attention that the question of whether the odd chromatic index of hypergraphs is bounded by a constant had already been posed in \cite{goetze2025strong}. This work was carried out as part of the first author’s PhD studies under the supervision of Professor Imre Leader, to whom the author is deeply indebted for his many insightful discussions and unwavering support.

G. Carenini is supported by the CB European PhD Studentship funded by Trinity College, Cambridge.

\bibliographystyle{unsrt}
\bibliography{references}

\end{document}